\documentclass[10pt]{article}
\usepackage{amsmath}
\usepackage{amssymb}
\usepackage{amsthm}
\usepackage{amsfonts}
\usepackage{algorithm}
\usepackage{algorithmic}

\DeclareMathOperator{\ord}{ord} \DeclareMathOperator{\qf}{qf}
\DeclareMathOperator{\Ker}{Ker} 
\DeclareMathOperator{\trdeg}{tr.deg}
\DeclareMathOperator{\Card}{Card}

\begin{document}

\def\D{\displaystyle}
\newtheorem{theorem}{Theorem}
\newtheorem{definition}{Definition}
\newtheorem{lemma}{Lemma}
\newtheorem{example}{Example}
\newtheorem{proposition}{Proposition}

\begin{center}
\large {\bf Difference Dimension Quasi-polynomials} \normalsize

\smallskip

Alexander  Levin

The Catholic University of America

Washington, D. C.  20064

levin@cua.edu

http://faculty.cua.edu/levin

\end{center}

\begin{abstract}
We consider Hilbert-type functions associated with difference (not
necessarily inversive) field extensions and systems of algebraic
difference equations in the case when the translations are assigned
some integer weights.  We will show that such functions are
quasi-polynomials, which can be represented as alternative sums of
Ehrhart quasi-polynomials associated with rational conic polytopes.
In particular, we obtain generalizations of main theorems on
difference dimension polynomials and their invariants to the case of
weighted basic difference operators.
\end{abstract}


\section{Introduction}
Difference dimension polynomials, first introduced in \cite{Levin1},
play the same role in the study of difference algebraic structures
and systems of algebraic difference equations as is played by
Hilbert polynomials in commutative algebra and algebraic geometry.
Most of the known results on such polynomials (including algorithms
for their computation) can be found in \cite{KLMP} and
\cite{Levin2}. Furthermore, as it is shown in \cite[Chapter
7]{Levin2}, difference dimension polynomials have an important
analytic interpretation: if a system of partial algebraic difference
equations represents a system of equations in finite differences,
then its difference dimension polynomial expresses the A. Einstein's
strength of the system. This fact determines the importance of the
study of difference dimension polynomials for the qualitative theory
of algebraic difference equations.

Another important application of difference dimension polynomials is
based on  the fact that such a polynomial can be naturally assigned
to any prime reflexive difference ideal of a finitely generated
difference algebra. The relationship between prime reflexive
difference ideals and their dimension polynomials has allowed one to
obtain new interesting results on the Krull-type dimension of
difference algebras and modules, as well as on local difference
algebras (see, for example, \cite{LM}, \cite[Section 4.6]{Levin2}
and \cite{Levin3}).

In this paper, we prove the existence and determine invariants of a
dimension quasi-polynomial associated with a difference field
extension with weighted basic translations. We also express a
difference dimension quasi-polynomial as an alternative sum of
Ehrhart quasi-polynomials associated with rational conic polytopes
and show that, given a system of algebraic difference equations with
weighted translations, one can assign to it a quasi-polynomial that
can be viewed as an algebraic version of the Einstein's strength of
the system. Note that systems of difference equations of these kind
arise, in particular, from finite difference approximations of
systems of PDEs with weighted derivatives, see, for example,
\cite{Shananin1} and \cite{Shananin2}. One should also mention that
the existence of Ehrhart-type dimension quasi-polynomials associated
with weighted filtrations of differential and inversive difference
modules was established by C. D\"onch in his dissertation
\cite{Donch}. C. D\"onch also proved the existence of dimension
quasi-polynomials associated with finitely generated differential
field extensions using the technique of weight Gr\"obner bases in
the associated modules of K\"ahler differentials. This approach,
however, cannot be applied in the case of non-inversive difference
fields; the main tool used in our paper is the method of
characteristic sets modified to the case of difference polynomials
over a difference field with weighted translations.

\section{Preliminaries}

Throughout the paper, $\mathbb{N}$, $\mathbb{Z}$, $\mathbb{Q}$, and
$\mathbb{R}$ denote the sets of all non-negative integers, integers,
rational numbers, and real numbers, respectively.  If $m$ is a
positive integer, then by the product order on $\mathbb{N}^{m}$ we
mean a partial order $<_{P}$ such that $(a_{1},\dots,
a_{m})<_{P}(a'_{1},\dots, a'_{m})$ if and only if $a_{i}<_{P}a'_{i}$
for $i=1,\dots, m$.

By a {\em difference ring} we mean a commutative ring $R$ together
with a finite set $\sigma = \{\alpha_{1},\dots, \alpha_{n}\}$ of
mutually commuting endomorphisms of $R$. The set $\sigma$ is called
a {\em basic set} of $R$ and the endomorphisms $\alpha_{i}$ are
called {\em translations}. We also say that $R$ is a $\sigma$-ring.
In what follows we assume that the translations are injective (this
assumption is standard in most works on difference algebra).
Furthermore, we will consider the free commutative semigroup of all
power products $\tau = \alpha_{1}^{k_{1}}\dots \alpha_{m}^{k_{m}}$
($k_{1},\dots, k_{m}\in \mathbb{N})$ denoted by $T$ (or
$T_{\sigma}$). If $\theta, \tau\in T$ and there is $\tau'\in T$ such
that $\tau = \theta\tau'$, we say that $\theta$ {\em divides} $\tau$
and write $\theta|\tau$. Otherwise, we write $\theta\nmid\tau$.

If a difference ($\sigma$-) ring is a field, it is called a
difference (or $\sigma$-) field.

A subring (ideal) $R_{0}$ of a $\sigma$-ring $R$ is said to be a
difference  (or $\sigma$-) subring of $R$ (respectively, difference
(or $\sigma$-) ideal of $R$) if $R_{0}$ is closed with respect to
the action of any translation $\alpha_{i}\in\sigma$. A
$\sigma$-ideal $I$ of a $\sigma$-ring $R$ is called {\em reflexive}
if the inclusion $\alpha_{i}(a)\in I$ ($a\in R,\, \alpha_{i}\in
\sigma$) implies the inclusion $a\in I$. In this case, the factor
ring $R/I$ has a natural structure of a difference ring with the
same basic set $\sigma$ where $\alpha(a+I) = \alpha(a) + I$ for any
coset $a+I\in R/I$ and $\alpha\in\sigma$.  If a difference ideal is
prime, it is referred to as a {\em prime difference ideal}.

If $R$ is a $\sigma$-ring and $S\subseteq R$, then the intersection
$I$ of all $\sigma$-ideals of $R$ containing the set $S$ is the
smallest $\sigma$-ideal of $R$ containing $S$; it is denoted by
$[S]$. Clearly, $[S]$ is generated, as an ideal, by the set
$\{\tau(a)\,|\, a\in S,\, \tau\in T\}$). If the set $S$ is finite,
$S = \{a_{1},\dots, a_{r}\}$, we say that the $\sigma$-ideal $I$ is
finitely generated (we write this as $I = [a_{1},\dots, a_{r}]$) and
call $a_{1},\dots, a_{r}$ difference (or $\sigma$-) generators of
$I$. If $R_{0}$ is a $\sigma$-subring of $R$ and $S\subseteq R$,
then the smallest $\sigma$-subring of $R$ containing $R_{0}$ and $S$
is denoted by $R_{0}\langle S \rangle$. As a ring, it is generated
over $R_{0}$ by the set $\{\tau s\,|\,\tau\in T,\, s\in S\}$. (Here
and below we frequently write$\tau s$ for $\tau(s)$.)

A ring homomorphism of $\sigma$-rings $\phi: R \longrightarrow S$ is
called a {\em difference} (or $\sigma$-) {\em homomorphism\/} if
$\phi(\alpha a) = \alpha \phi(a)$ for any $\alpha\in\sigma$, $a\in
R$. It is easy to see that the kernel of such a mapping is a
reflexive difference ideal of $R$.

If $L$ is a difference ($\sigma$-) field and its subfield $K$ is
also a $\sigma$-subring of $L$, then  $K$ is said to be a difference
(or $\sigma$-) subfield of $L$; $L$, in turn, is called a difference
(or $\sigma$-) field extension or a $\sigma$-overfield of $K$. We
also say that we have a $\sigma$-field extension $L/K$.

If $K$ is a $\sigma$-subfield of a $\sigma$-field $L$ and
$S\subseteq L$, then the smallest $\sigma$-subfield of $L$
containing $K$ and $S$ is denoted by $K\langle S \rangle$. If the
set $S$ is finite, $S = \{\eta_{1},\dots,\eta_{n}\}$, then $K\langle
S \rangle$ is written as $K\langle \eta_{1},\dots,\eta_{n}\rangle$
and is said to be a finitely generated difference (or $\sigma$-)
extension of $K$ with the set of $\sigma$-generators
$\{\eta_{1},\dots,\eta_{n}\}$. It is easy to see that the field
$K\langle \eta_{1},\dots,\eta_{n}\rangle$ coincides with the field
$K_0(\{\tau\eta_{i}\,|\,\tau\in T,\, 1\leq i\leq n\}$).

If $R$ is a difference ($\sigma$-) ring and $Y =\{y_{1},\dots,
y_{n}\}$ is a finite set of symbols, then one can consider the
countable set of symbols $TY = \{\tau y_{j}|\tau\in T,  1\leq j\leq
n\}$ and the polynomial ring $R[TY]$ in the set of indeterminates
$TY$ over $R$. This polynomial ring is naturally viewed as a
$\sigma$-ring where the action of the translations on $R$ is
extended to $R[TY]$ by setting $\alpha(\tau y_{j}) =
(\alpha\tau)y_{j}$ for any $\alpha\in\sigma$, $\tau\in T$, $1\leq
j\leq n$. The constructed difference ring is denoted by
$R\{y_{1},\dots, y_{n}\}$; it is called the {\em ring of difference}
(or $\sigma$-) {\em polynomials} in the set of difference
($\sigma$-) indeterminates $y_{1},\dots, y_{n}$ over $R$. Elements
of $R\{y_{1},\dots, y_{n}\}$ are called difference (or $\sigma$-)
polynomials. If $R$ is a $\sigma$-subring of a $\sigma$-ring $S$,
$A\in R\{y_{1},\dots, y_{n}\}$ and $\eta = (\eta_{1},\dots,
\eta_{n})\in S^{n}$, then $A(\eta)$ denotes the result of the
replacement of every entry $\tau y_{i}$ in $A$ by $\tau \eta_{i}$
($\tau\in T$, $1\leq i\leq n$).

Let $U = \big\{u^{(\lambda)} | \lambda \in \Lambda\big\}$ be a
family of elements in some $\sigma$-overring of a difference
($\sigma$-) ring $R$. We say that the family $U$ is {\em
transformally} (or $\sigma$-{\em algebraically) dependent} over $R$,
if the family
$$U^{\sigma} = \left\{\tau u^{(\lambda)}\,|\,\tau\in T,\, \lambda \in \Lambda\right\}$$
is algebraically dependent over $R$ (that is, there exist elements
$v^{(1)},\dots, v^{(k)}\in U^{\sigma}$ and a nonzero polynomial $f$
in $k$ variables with coefficients from $R$ such that
$f\big(v^{(1)},\dots, v^{(k)}\big) = 0$). Otherwise, the family $U$
is said to be {\em transformally} (or $\sigma$-{\em algebraically)
independent} over $R$.

If $K$ is a difference ($\sigma$-) field and $L$ a $\sigma$-field
extension of $K$, then a set $B\subseteq L$ is said to be a {\em
difference} (or $\sigma$-) {\em transcendence basis} of $L$ over $K$
if $B$ is $\sigma$-algebraically independent over $K$ and every
element $a\in L$ is $\sigma$-algebraic over $K\langle B\rangle$
(that is, the set $\{\tau a\,|\,\tau\in T\}$ is algebraically
dependent over $K$) . If $L$ is a finitely generated $\sigma$-field
extension of $K$, then all $\sigma$-transcendence bases of $L$ over
$K$ are finite and have the same number of elements (see
\cite[Section 4.1]{Levin2}). This number is called the {\em
difference} (or $\sigma$-) {\em transcendence degree} of $L$ over
$K$ (or the $\sigma$-transcendence degree of the extension $L/K$);
it is denoted by $\sigma$-$\trdeg_{K}L$.

\section{Dimension quasi-polynomials of subsets of
$\mathbb{N}^{m}$}

A function $f:\mathbb{Z}\rightarrow \mathbb{Q}$ is called a
(univariate) {\em quasi-polynomial} of period $q$ if there exist $q$
polynomials $g_{i}(x)\in\mathbb{Q}[x]$ \,\,($0\leq i\leq q-1$) such
that $f(n) = g_{i}(n)$ whenever $n\in\mathbb{Z}$ and $n\equiv q\,
(mod\, n)$.

An equivalent way of introducing quasi-polynomials is as follows.

A {\em rational periodic number} $U(n)$ is a function
$U:\mathbb{Z}\rightarrow \mathbb{Q}$ with the property that there
exists (a period) $q\in\mathbb{N}$ such that
$$U(n) = U(n')\,\,\, {\text{whenever}}\,\,\,
n\equiv n' \,(mod\, q).$$

A rational periodic number can be represented by a list of $q$ its
possible values enclosed in square brackets: $$U(n) = [a_{0},\dots,
a_{q}]_{n}.$$ For example, $U(n) = \left[{\D\frac{1}{2}},\,
{\D\frac{3}{4}},\, 1\right]_{n}$ is a periodic number with period
$3$ such that $U(n) = {\D\frac{1}{2}}$ if $n\equiv 0 \,(mod\, 3)$,
$U(n) = {\D\frac{3}{4}}$ if $n\equiv 1 \,(mod\, 3)$, and  $U(n) = 1$
if $n\equiv 2 \,(mod\, 3)$.

\medskip

With the above notation, a (univariate) {\em quasi-polynomial of
degree $d$} is a function $f:\mathbb{Z}\rightarrow \mathbb{Q}$ such
that
$$f(n) = c_{d}(n)n^{d} + \dots + c_{1}(n)n + c_{0}(n)\,\,\,\,\,
(n\in \mathbb{Z})$$ where $c_{i}(n)$'s are rational periodic numbers
and $c_{d}(n)\neq 0$.

\medskip

One of the main applications of the theory of quasi-polynomials is
its application to the problem of counting integer points in
polytopes.

Recall that a {\em rational polytope} in $\mathbb{R}^{d}$ is the
convex hall of finitely many points (vertices) in $\mathbb{Q}^{d}$.
Equivalently, a rational polytope $P\subseteq \mathbb{R}^{d}$ is the
set of solutions of a finite system of linear inequalities $$A{\bf
x}\leq {\bf b},$$ where $A$ is an $m\times d$-matrix with integer
entries ($m$ is a positive integer) and ${\bf b}\in \mathbb{Z}^{m}$,
provided that the solution set is bounded. If all vertices of a
rational polytope have integer coordinates, it is called {\em
lattice}.

\medskip

Let $P\subseteq\mathbb{R}^{d}$ be a rational polytope. (We assume
that $P$ has dimension $d$, that is, $P$ is not contained in a
proper affine subspace of $\mathbb{R}^{d}$.) Then a polytope
$$rP = \{r{\bf x}\,|\,{\bf x}\in P\}$$ ($r\in \mathbb{N}$, $n\geq
1$) is called the  $r$th {\em dilate} of $P$. Clearly, if ${\bf
v}_{1},\dots, {\bf v}_{k}$ are all vertices of $P$, then $rP$ is the
convex hall of $r{\bf v}_{1},\dots, r{\bf v}_{k}$.

Given a rational polytope $P$, let, $L(P, r)$ denote the number of
integer points (that is, points with integer coordinates) in $rP$
(in other words, $L(P, r) = \Card(rP\cap \mathbb{Z}^{d}$). The
following result is due to E. Ehrhart, see \cite{Ehrhart}.

\begin{theorem}
Let $P\subseteq\mathbb{R}^{d}$ be a rational polytope. Then

\smallskip

{\em (i)} \,$L(P, r)$ is a degree $d$ quasi-polynomial of $r$.

\smallskip

{\em (ii)} \,The coefficient of the leading term of this
quasi-polynomial is equal to the Euclidean volume of $P$.

\smallskip

{\em (iii)}\, The minimum period of $L(P, r)$ is a divisor of the
number $\mathcal{D}(P)= \min\{n\in\mathbb{N}\,|\,nP$ is a lattice
polytope$\}$.

\smallskip

(iv) \,If $P$ is a lattice polynomial, then $L(P, r)$ is a
polynomial of $n$ with rational coefficients.
\end{theorem}
The main tools for computation of Ehrhart quasi-polynomials are
Alexander Barvinok's polynomial time algorithm and its
modifications, see \cite{Barvinok1},  \cite{Barvinok2} and
\cite{Barvinok3}. In some cases, however, the Ehrhart
quasi-polynomial can be found directly from the Ehrhart's theorem by
evaluating the periodic numbers, which are coefficients of the
quasi-polynomial.

\begin{example}
{\em Consider a polytope $$P = \{(x_{1}, x_{2})\in
\mathbb{R}^{2}\,|\,x_{1}\geq 0,\, x_{2}\geq 0,\, x_{1}+x_{2}\leq
3,\, 2x_{1}\leq 5\}.$$

Then $$rP = \{(x_{1}, x_{2})\in \mathbb{R}^{2}\,|\,x_{1}\geq 0,\,
x_{2}\geq 0,\, x_{1}+x_{2}\leq 3r,\, 2x_{1}\leq 5r\}$$ is a polytope
with vertices $(0, 0)$, $\left({\D\frac{5r}{2}}, 0\right)$, $(0,
3r)$, and $\left({\D\frac{5r}{2}}, {\D\frac{r}{2}}\right)$. By the
Ehrhart Theorem, $$L(P, r) = \alpha r^{2} +
[\beta_{1},\beta_{2}]_{r}r + [\gamma_{1},\gamma_{2}]_{r}.$$

The direct counting gives $L(P, 0) = 1$, $L(P, 1) = 9$, $L(P, 2) =
27$, $L(P,3) = 52$, and $L(P, 4) = 88$. The following figure shows
integer points in $rP$ for $r = 0, 1, 2$.

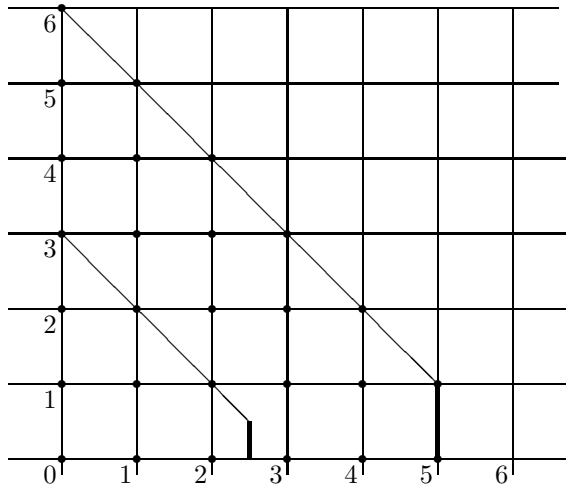
\setlength{\unitlength}{1cm}
\begin{picture}(12,6.5)
\put(1.6,0.3){\line(1,0){0.7}} \put (2.3,0.30){\circle*{0.1}}
\put(2.3,0.3){\line(1,0){1}} \put (3.3,0.30){\circle*{0.1}}
\put(4.3,0.30){\circle*{0.1}} \put (6.3,0.30){\circle*{0.1}}
\put(5.3,0.30){\circle*{0.1}} \put(3.3,0.3){\line(1,0){1}}
\put(4.3,0.3){\line(1,0){1}} \put(5.3,0.3){\line(1,0){1}}
\put(6.3,0.3){\line(1,0){1}} \put (7.3,0.30){\circle*{0.1}}
\put(7.3,0.3){\line(1,0){1}} \put(8.3,0.3){\line(1,0){0.7}}
\put(2.3,0.3){\line(0,-1){0.2}} \put(3.3,0.3){\line(0,-1){0.2}}
\put(4.3,0.3){\line(0,-1){0.2}} \put(5.3,0.3){\line(0,-1){0.2}}
\put(6.3,0.3){\line(0,-1){0.2}} \put(7.3,0.3){\line(0,-1){0.2}}
\put(8.3,0.3){\line(0,-1){0.2}} \put(2.3,0.3){\line(0,1){1}}
\put(3.3,0.3){\line(0,1){1.5}} \put(4.3,0.3){\line(0,1){1}}
\put(5.3,0.3){\line(0,1){1}} \put(6.3,0.3){\line(0,1){1}}
\put(7.3,0.3){\line(0,1){1}} \put(8.3,0.3){\line(0,1){1}}
\put(7.3,1.30){\circle*{0.1}} \put(2.3,1.30){\circle*{0.1}}
\put(3.3,1.30){\circle*{0.1}} \put(4.3,1.30){\circle*{0.1}}
\put(5.3,1.30){\circle*{0.1}} \put(6.3,1.30){\circle*{0.1}}
\put(1.6,1.3){\line(1,0){0.7}} \put(2.3,1.3){\line(1,0){1}}
\put(3.3,1.3){\line(1,0){1}} \put(4.3,1.3){\line(1,0){1}}
\put(5.3,1.3){\line(1,0){1}} \put(6.3,1.3){\line(1,0){1}}
\put(7.3,1.3){\line(1,0){1}} \put(8.3,1.3){\line(1,0){0.7}}
\put(2.3,2.3){\circle*{0.1}} \put(3.3,2.3){\circle*{0.1}}
\put(4.3,2.3){\circle*{0.1}} \put(5.3,2.3){\circle*{0.1}}
\put(6.3,2.3){\circle*{0.1}} \put(2.3,3.3){\circle*{0.1}}
\put(3.3,3.3){\circle*{0.1}} \put(4.3,3.3){\circle*{0.1}}
\put(5.3,3.3){\circle*{0.1}} \put(2.3,5.3){\circle*{0.1}}
\put(2.3,6.3){\circle*{0.1}} \put(2.3,4.3){\circle*{0.1}}
\put(3.3,4.3){\circle*{0.1}} \put(4.3,4.3){\circle*{0.1}}
\put(1.6,2.3){\line(1,0){0.7}} \put(2.3,2.3){\line(1,0){1}}
\put(3.3,2.3){\line(1,0){1}} \put(4.3,2.3){\line(1,0){1}}
\put(5.3,2.3){\line(1,0){1}} \put(6.3,2.3){\line(1,0){1}}
\put(7.3,2.3){\line(1,0){1}} \put(8.3,2.3){\line(1,0){0.7}}
\put(1.6,3.3){\line(1,0){0.7}} \put(2.3,3.3){\line(1,0){1}}
\put(3.3,3.3){\line(1,0){1}} \put(4.3,3.3){\line(1,0){1}}
\put(5.3,3.3){\line(1,0){1}} \put(6.3,3.3){\line(1,0){1}}
\put(7.3,3.3){\line(1,0){1}} \put(8.3,3.3){\line(1,0){0.7}}
\put(1.6,4.3){\line(1,0){0.7}} \put(2.3,4.3){\line(1,0){1}}
\put(3.3,4.3){\line(1,0){1}} \put(4.3,4.3){\line(1,0){1}}
\put(5.3,4.3){\line(1,0){1}} \put(6.3,4.3){\line(1,0){1}}
\put(7.3,4.3){\line(1,0){1}} \put(8.3,4.3){\line(1,0){0.7}}
\put(2.3,1.3){\line(0,1){1}} \put(3.3,1.3){\line(0,1){1}}
\put(4.3,1.3){\line(0,1){1}} \put(5.3,1.3){\line(0,1){1}}
\put(6.3,1.3){\line(0,1){1}} \put(7.3,1.3){\line(0,1){1}}
\put(8.3,1.3){\line(0,1){1}} \put(2.3,2.3){\line(0,1){1}}
\put(3.3,2.3){\line(0,1){1}} \put(4.3,2.3){\line(0,1){1}}
\put(5.3,2.3){\line(0,1){1}} \put(6.3,2.3){\line(0,1){1}}
\put(7.3,2.3){\line(0,1){1}} \put(8.3,2.3){\line(0,1){1}}
\put(2.3,3.3){\line(0,1){1}} \put(3.3,3.3){\line(0,1){1}}
\put(4.3,3.3){\line(0,1){1}} \put(5.3,3.3){\line(0,1){1}}
\put(6.3,3.3){\line(0,1){1}} \put(7.3,3.3){\line(0,1){1}}
\put(8.3,3.3){\line(0,1){1}} \put(2.3,4.3){\line(0,1){0.5}}
\put(3.3,4.3){\line(0,1){0.5}} \put(4.3,4.3){\line(0,1){0.5}}
\put(5.3,4.3){\line(0,1){0.5}} \put(6.3,4.3){\line(0,1){0.5}}
\put(7.3,4.3){\line(0,1){0.5}} \put(8.3,4.3){\line(0,1){0.5}}
\put(2.15,2.1){\makebox(0,0){$2$}}
\put(2.15,1.1){\makebox(0,0){$1$}}
\put(5.15,0.1){\makebox(0,0){$3$}}
\put(4.15,0.1){\makebox(0,0){$2$}}
\put(3.15,0.1){\makebox(0,0){$1$}}
\put(2.15,0.1){\makebox(0,0){$0$}}
\put(6.15,0.1){\makebox(0,0){$4$}}
\put(7.15,0.1){\makebox(0,0){$5$}}
\put(8.15,0.1){\makebox(0,0){$6$}}
\put(2.15,4.1){\makebox(0,0){$4$}}
\put(2.15,5.1){\makebox(0,0){$5$}}
\put(2.15,6.1){\makebox(0,0){$6$}}
\put(2.15,3.1){\makebox(0,0){$3$}} \put(1.6,5.3){\line(1,0){7.3}}
\put(1.6,6.3){\line(1,0){7.3}} \put(3.3,5.3){\circle*{0.1}}
\put(2.3,4.3){\line(0,1){1}} \put(2.3,5.3){\line(0,1){1}}
\put(3.3,4.3){\line(0,1){1}} \put(4.3,4.3){\line(0,1){1}}
\put(5.3,4.3){\line(0,1){1}} \put(6.3,4.3){\line(0,1){1}}
\put(7.3,4.3){\line(0,1){1}} \put(8.3,4.3){\line(0,1){1}}
\put(3.3,5.3){\line(0,1){1}} \put(4.3,5.3){\line(0,1){1}}
\put(5.3,5.3){\line(0,1){1}} \put(6.3,5.3){\line(0,1){1}}
\put(7.3,5.3){\line(0,1){1}} \put(8.3,5.3){\line(0,1){1}}

\linethickness{1.5pt} \put(2.3,6.3){\line(1,-1){5}}
\put(2.3,3.3){\line(1,-1){2.5}} \put(4.8,0.3){\line(0,1){0.5}}
\put(7.3,0.3){\line(0,1){1}} \thinlines
\end{picture}

Substituting $r=0, 1, 2, 3, 4$ into the expression $$L(P, r) =
\alpha r^{2} + [\beta_{1},\beta_{2}]_{r}r +
[\gamma_{1},\gamma_{2}]_{r}$$ one obtains a system of linear
equations
\[ \begin{cases}
\gamma_{1} = 1,
\\ \alpha + \beta_{2} + \gamma_{2} = 9,\\
 4\alpha + 2\beta_{1} + \gamma_{1} = 27,\\
9\alpha + 3\beta_{2} + \gamma_{2} = 52,\\
16\alpha + 4\beta_{1} + \gamma_{1} = 88.
\end{cases}\]

That gives $\alpha = {\D\frac{38}{5}}$, $\beta_{1} =
{\D\frac{17}{4}}$, $\beta_{2} = 4$, $\gamma_{1} = 1$, and
$\gamma_{2} = {\D\frac{5}{8}}$. Thus,
$$L(P, r) = {\D\frac{38}{5}}r^{2} + \left[{\D\frac{17}{4}}, 4\right]_{r}r + \left[1,
{\D\frac{5}{8}}\right]_{r}.$$}
\end{example}

\bigskip

Let $w_{1},\dots, w_{m}$ be fixed positive integers ($m > 0$). Then
for any $m$-tuple $a = (a_{1},\dots , a_{m})\in \mathbb{N}^{m}$, the
number $$\ord_{w}a = w_{1}a_{1} + \dots + w_{m}a_{m}$$ is called the
{\em order of $a$ with respect to the weights} $w_{1},\dots, w_{m}$.
If the weights are fixed, $\ord_{w}a$ is simply called the {\em
order of $a$}.

In what follows, $\lambda^{(m)}_{w}(t)$ denotes the Ehrhart
quasi-polynomial that describes the number of integer points in the
{\em conic polytope} defined as $\{(x_{1},\dots, x_{m})\in
\mathbb{R}^{m}\,|\, \sum_{i=1}^{m}w_{i}x_{i}\leq t,\, x_{j}\geq 0\,
(1\leq j\leq m)\}$. It follows from the Ehrhart's Theorem that
$\lambda^{(m)}_{w}(t)$ is a quasi-polynomial of degree $m$ whose
leading coefficient is $\D\frac{1}{m!w_{1}\dots w_{m}}$. A
polynomial time algorithm for computing $\lambda^{(m)}_{w}(t)$ can
be found, for example, in \cite{Barvinok3}.

If $A\subseteq \mathbb{N}^{m}$ and $r\in \mathbb{N}$, we set
$$A^{(w)}(r) = \{a = (a_{1},\dots, a_{m})\in A\,|\,\ord_{w}a \leq r\}.$$
Furthermore, if $A\subseteq \mathbb{N}^{m}$,
then $V_{A}$ will denote the set of all $m$-tuples $v = (v_{1},\dots
, v_{m})\in\mathbb{N}^{m}$ that are not greater than or equal to any
$m$-tuple from $A$ with respect to $\leq_{P}$. (Recall that
$\leq_{P}$ denotes the product order on $\mathbb{N}^{m}$.)

In other words, an element $v=(v_{1}, \dots , v_{m})\in
\mathbb{N}^{m}$ belongs to $V_{A}$ if and only if for any element
$(a_{1},\dots , a_{m})\in A$ there exists $i\in \mathbb{N}, 1\leq
i\leq m$, such that $a_{i} > v_{i}$.

The following theorem generalizes E. Kolchin's result on dimension
polynomials of subsets of $\mathbb{N}^{m}$ (see \cite[Chapter 0,
Lemma 16]{Kolchin}

\begin{theorem}  With the above conventions, for any set
$A\subseteq \mathbb{N}^{m}$, there exists a quasi-polynomial
$\chi^{(w)}_{A}(t)$ in one variable $t$ such that

\smallskip

{\em (i)} \, $\chi^{(w)}_{A}(r) = \Card\,V^{(w)}_{A}(r)$ for all
sufficiently large $r\in\mathbb{N}$.

\smallskip

{\em (ii)} \, $\deg\,\chi^{(w)}_{A}(t)\leq m$.

\smallskip

{\em (iii)} \, $\deg\,\chi^{(w)}_{A}(t) = m$ if and only if $A =
\emptyset$. In this case $\chi^{(w)}_{A}(t) = \lambda^{(m)}_{w}(t)$.

\smallskip

{\em (iv)}  \, $\chi^{(w)}_{A}(t) = 0$ if and only if $(0,\dots,
0)\in A$.
\end{theorem}

PROOF.\, Let $A\subseteq \mathbb{N}^{m}$ and $V_{A} = \{v\in
\mathbb{N}^{m}\,|\,a\nleq_{P}v$ for any $a\in A\}$. Clearly, if one
replaces $A$ with the finite set of all minimal with respect to
$\leq_P$ points of $A$, this replacement does not change $V_{A}$.
Therefore, we can assume that $A$ is finite. Let $A = \{a^{(1)},
\dots, a^{(d)}\}$ where $a^{(i)} = (a_{i1},\dots, a_{im})$.

We proceed by induction on $m$ (for $m =1$ the statement is
obviously true) and $|A| = \D\sum_{i=1}^{d}\D\sum_{j=1}^{m}a_{ij}$.

For any $s\in \mathbb{N}$ and $B\subseteq \mathbb{N}^{m}$, let
$$N_{B}(s) = \Card\,V^{(w)}_{B}(s) = \Card\{v = (v_{1},\dots, v_{m})\in V_{B}\,|\,\ord_{w}v
\leq s\}.$$

If $v = (v_{1},\dots, v_{m})\in V^{(w)}_{A}(s)$, then either

(a)\, $v_{m} = 0$\, and \, $\sum_{i=1}^{m}w_{i}v_{i}\leq s$\, or

(b)\,  $v_{m} = v'_{m}+1$\, with\, $v'_{m}\in \mathbb{N}$\, and\,
$\sum_{i=1}^{m-1}w_{i}v_{i} + w_{m}v'_{m}\leq s-w_{m}$.

Let
$$A_{0} = \{a = (a_{1},\dots, a_{m-1})\in \mathbb{N}^{m-1}\,|\,(a_{1},\dots, a_{m-1},
0)\in A\}\hspace{0.1in}\text{and}$$

$A_{1} = \{a = (a_{1},\dots, a_{m-1}, a'_{m})\in \mathbb{N}^{m}\,|\,
(a_{1},\dots, a_{m-1}, a'_{m}+1)\in A$ or

$a'_{m} = 0$ and $(a_{1},\dots, a_{m-1}, 0)\in A\}.$

\medskip

It is easy to see that there is a one-to-one correspondence between
elements of $V^{(w)}_{A}(s)$ satisfying condition (a) and elements
of the set $V^{(w)}_{A_{0}}(s)\subseteq \mathbb{N}^{m-1}$
(\,$(v_{1},\dots, v_{m-1}, 0) \longleftrightarrow (v_{1},\dots,
v_{m-1})$\,). Also, there is a one-to-one correspondence between
elements of $V^{(w)}_{A}(s)$ satisfying condition (b) and elements
of the set $V^{(w)}_{A_{1}}(s-w_{m})$ (\,$(v_{1},\dots, v_{m-1},
v'_{m}+1) \longleftrightarrow (v_{1},\dots, v_{m-1}, v'_{m})$\,). It
follows that $$N_{A}(s) = N_{A_{0}}(s) + N_{A_{1}}(s-w_{m}).$$

Since $|A_{1}| < |A|$ and $A_{0}\subseteq\mathbb{N}^{m-1}$,
$N_{A_{0}}(s)$ and $N_{A_{1}}(s-w_{m})$ are expressed by
quasi-polynomials of degree at most $m$. Therefore, $N_{A}(s)$ is
expressed by such a quasi-polynomial as well.

If $A = \emptyset$, then $V_{A} = \mathbb{N}^{m}$ and $\omega_{A}(t)
= \lambda^{(m)}_{w}(t)$. In order to complete the proof of part
(iii) of the theorem, we need to show that if $A\neq\emptyset$, $A =
\{a^{(1)}, \dots, a^{(d)}\}$ with $a^{(i)} = (a_{i1},\dots,
a_{im})$, then $\deg\,\chi^{(w)}_{A} < m$.

For every $j\in\{1,\dots, m\}$, let $e_{j} = \min\{a_{ij}\,|\,1\leq
i\leq d\}$ and let $e = (e_{1},\dots, e_{m})$. Then $V_{A}\subseteq
V_{\{e\}}$, hence $\chi^{(w)}_{A}(r)\leq \chi^{(w)}_{\{e\}}(r)$ for
all sufficiently large $r\in \mathbb{N}$. Therefore, it remains to
notice that $\deg\,\chi^{(w)}_{\{e\}} < m$. Indeed, for all
sufficiently large $r\in \mathbb{N}$, $$V_{\{e\}}(r) =
V_{\emptyset}(r)\setminus\{b\in \mathbb{N}^{m}\,|\,e\leq_{P}b\,\,
\text{and}\,\, \ord_{w}b\leq r\}.$$ Therefore,
$$\chi^{(w)}_{\{e\}}(t) = \lambda^{(m)}_{w}(t) -
\lambda^{(m)}_{w}(t-\ord_{w}e).$$ Since both quasi-polynomials
$\lambda^{(m)}_{w}(t)$ and $\lambda^{(m)}_{w}(t-\ord_{w}e)$ are of
degree $m$ and have the same leading coefficient
$\D\frac{1}{m!w_{1}\dots w_{m}}$, $\deg\,\chi^{(w)}_{\{e\}} < m$.

The last part of the theorem follows from the observation that
$(0,\dots, 0)\in A$ if and only if $V_{A}=\emptyset$, which is
equivalent to the equality $\chi^{(w)}_{A}(t) = 0$.\, $\square$

\begin{definition} The quasi-polynomial $\chi^{(w)}_{A}(t)$ whose
existence is established by Theorem 2 is called the dimension
quasi-polynomial of the set $A\subseteq \mathbb{N}^{m}$ associated
with the weight vector $(w_{1},\dots, w_{m})$.
\end{definition}

As we have mentioned, the dimension quasi-polynomial of a set
$A\subseteq \mathbb{N}^{m}$ coincides with the dimension
quasi-polynomial of the finite set of all minimal points of $A$ with
respect to the product order on $\mathbb{N}^{m}$. The following
theorem gives a formula for computation of the dimension
quasi-polynomial of a finite set. It generalizes the corresponding
result for dimension polynomials of finite subsets of
$\mathbb{N}^{m}$, see \cite[Proposition 2.3.1]{MP}.

\begin{theorem}
Let $A = \{a^{(1)} ,\dots, a^{(d)}\}$ be a finite subset of
$\mathbb{N}^{m}$ and let $a^{(i)} = (a_{i1},\dots, a_{im})$ for
$i=1,\dots, d$. For any $l = 0,\dots, d$, let $\Gamma(l, d)$ denote
the set of all $l$-element subsets of $\{1,\dots, d\}$ ($\Gamma(0,
d) = \emptyset$) and for any set $\epsilon = \{a^{(i_{1})},\dots,
a^{(i_{p})}\}\in \Gamma(l, d)$ ($1\leq i_{1} <\dots < i_{p}\leq d$),
let $c_{\epsilon j} = \max\{a_{\nu j}\,|\,1\leq \nu\leq p\}$ (the
maximal $j$th coordinate of elements of $\epsilon$). Furthermore,
let $c_{\epsilon} = (c_{\epsilon 1}, \dots, c_{\epsilon m})$. Then
\begin{equation}
\chi^{(w)}_{A}(t) = \sum_{l=0}^{d}(-1)^{l}\sum_{\epsilon\in\Gamma(l,
d)}\lambda^{(m)}_{w}(t - \ord_{w}c_{\epsilon}).
\end{equation}
\end{theorem}

PROOF.\, Let $\epsilon\in\Gamma(l, d)$ and for any $s\in\mathbb{N}$,
let $$C_{\epsilon}(s) = \{c = (c_{1},\dots,
c_{m})\in\mathbb{N}^{m}\,|\,c_{\epsilon}\leq_{P} c\,\,
\text{and}\,\, \ord_{w}c\leq s\}.$$ Taking into account that
$\Card\,C_{\epsilon}(s) = \lambda^{(m)}_{w}(t -
\ord_{w}c_{\epsilon})$ for all $s\geq\ord_{w}c_{\epsilon}$, we can
mimic the proof of Proposition 2.3.1 of \cite{MP} as follows.

Let us show that the sum in the right-hand side of (1) counts all
points of $V^{(w)}_{A}(s)$ and every point of this set is counted
exactly once. Indeed, let a point $b = (b_{1},\dots, b_{m})\in
\mathbb{N}^{m}(s)$ exceed exactly $k$ points from the set $A$ ($k$
can take values from $0$ to $d$). If $k=0$, then $b\in
V^{(w)}_{A}(s)$ and $b$ is taken into account exactly once in the
term $(-1)^{0}\lambda^{(m)}_{w}(s)$ in (1).

Let $k > 0$. If we fix $l$ ($0\leq l\leq d$), then the point $b$
falls into exactly ${k\choose l}$ subsets of the form
$C_{\epsilon}(s)$ with $\epsilon\in\Gamma(l, d)$ (we set ${k\choose
l} = 0$ if $l> k$), so in the sum (1) the point $b$ is taken into
account $(-1)^{l}{k\choose l}$ times in the term
$(-1)^{l}\lambda^{(m)}_{w}(t - \ord_{w}c_{\epsilon})$. Therefore,
the total number of times the point $b$ is counted in the right-hand
side of (1) is $\sum_{l=0}^{k}(-1)^{l}{k\choose l} = (1-1)^{k} = 0$.
Thus, for all sufficiently large $s\in \mathbb{N}$, both right-hand
and left-hand sides of (1) give the number of points in the set
$V^{(w)}_{A}(s)$. This completes the proof of the theorem.\,
$\square$

The following example illustrates the above approach to the
computation of a dimension polynomial in the case of a subset of
$\mathbb{N}^{2}$.

\begin{example} {\em   Let $\mathcal{A} = \{((2,1), (0, 3)\}\subset
\mathbb{N}^{2}$ and let $w_{1} = 2,\, w_{2} = 1.$ Then the set
$V_{\mathcal{A}}$ consists of all integer points that lie in the
shadowed region in the following figure, where $C(2, r-4)$,
$D({\frac{r-3}{2}}, 3)$, and $E({\frac{r-1}{2}}, 1)$. ($BF$ is a
segment of the line $2x + y = r$.)

\setlength{\unitlength}{0.6cm}
\begin{picture}(12,11)

 \put(3.3,0.3){\line(1,0){1}}
\put(4.3,0.3){\line(1,0){1}} \put(5.3,0.3){\line(1,0){1}}
\put(6.3,0.3){\line(1,0){1}} \put(7.3,0.3){\line(1,0){1}}
\put(8.3,0.3){\line(1,0){0.7}} \put(2.3,0.3){\line(0,-1){0.2}}
\put(3.3,0.3){\line(0,-1){0.2}} \put(4.3,0.3){\line(0,-1){0.2}}
\put(5.3,0.3){\line(0,-1){0.2}} 
\put(7.3,0.3){\line(0,-1){0.2}} \put(8.3,0.3){\line(0,-1){0.2}}
\put(2.3,0.3){\line(0,1){1}} \put(3.3,0.3){\line(0,1){1.5}}
\put(2.3,7.3){\line(0,1){1}}\put(2.3,6.3){\line(0,1){1}}
\put(1.6,0.3){\line(1,0){0.7}} \put(2.3,0.3){\line(1,0){1}}
\put(4.3,0.3){\line(0,1){1}} \put(5.3,0.3){\line(0,1){1}}
\put(7.3,0.3){\line(0,1){1}} \put(8.3,0.3){\line(0,1){1}}
\put(1.6,1.3){\line(1,0){0.7}} \put(2.3,1.3){\line(1,0){1}}
\put(3.3,1.3){\line(1,0){1}} \put(4.3,1.3){\line(1,0){1}}
\put(5.3,1.3){\line(1,0){1}} \put(6.3,1.3){\line(1,0){1}}
\put(7.3,1.3){\line(1,0){1}} \put(8.3,1.3){\line(1,0){0.7}}
\put(1.6,8.3){\line(1,0){7.4}} \put(1.6,9.3){\line(1,0){7.4}}
\put(1.6,7.3){\line(1,0){7.4}} \put(1.6,2.3){\line(1,0){0.7}}
\put(2.3,2.3){\line(1,0){1}} \put(3.3,2.3){\line(1,0){1}}
\put(4.3,2.3){\line(1,0){1}} \put(5.3,2.3){\line(1,0){1}}
\put(6.3,2.3){\line(1,0){1}} \put(7.3,2.3){\line(1,0){1}}
\put(8.3,2.3){\line(1,0){0.7}} \put(1.6,3.3){\line(1,0){0.7}}

\linethickness{1pt} \put(2.3,10.3){\line(1,-2){5}}
\put(2.3,3.3){\line(1,0){1}} \put(3.3,3.3){\line(1,0){2.5}}
\put(4.3,1.3){\line(0,1){2}} \put(4.3,1.3){\line(1,0){2.5}}
\put(2.3,0.3){\line(0,1){10}} \put(2.3,0.3){\line(1,0){5}}
\put(4.3,3.3){\line(0,1){3}}

\thinlines \put(4.3,3.3){\line(1,0){1}}

\put(5.3,3.3){\line(1,0){1}} \put(6.3,3.3){\line(1,0){1}}
\put(7.3,3.3){\line(1,0){1}} \put(8.3,3.3){\line(1,0){0.7}}
\put(1.6,4.3){\line(1,0){0.7}} \put(2.3,4.3){\line(1,0){1}}
\put(3.3,4.3){\line(1,0){1}} \put(4.3,4.3){\line(1,0){1}}
\put(5.3,4.3){\line(1,0){1}} \put(6.3,4.3){\line(1,0){1}}
\put(7.3,4.3){\line(1,0){1}} \put(8.3,4.3){\line(1,0){0.7}}
\put(2.3,1.3){\line(0,1){1}} \put(3.3,1.3){\line(0,1){1}}
 \put(5.3,1.3){\line(0,1){9}} \put(7.3,1.3){\line(0,1){1}}
 \put(8.3,1.3){\line(0,1){1}} \put(2.3,2.3){\line(0,1){1}}
 \put(3.3,2.3){\line(0,1){1}}
\put(4.3,2.3){\line(0,1){1}} \put(5.3,2.3){\line(0,1){1}}
\put(7.3,2.3){\line(0,1){1}} \put(8.3,2.3){\line(0,1){1}}
\put(2.3,3.3){\line(0,1){1}} \put(3.3,3.3){\line(0,1){1}}
\put(4.3,3.3){\line(0,1){1}} \put(5.3,3.3){\line(0,1){1}}
\put(7.3,3.3){\line(0,1){1}} \put(8.3,3.3){\line(0,1){1}}
\put(2.3,4.3){\line(0,1){0.5}} \put(3.3,4.3){\line(0,1){0.5}}
\put(4.3,4.3){\line(0,1){0.5}} \put(5.3,4.3){\line(0,1){0.5}}
\put(7.3,4.3){\line(0,1){0.5}} \put(8.3,4.3){\line(0,1){0.5}}
\put(1.6,5.3){\line(1,0){7.3}} \put(1.6,6.3){\line(1,0){7.3}}
\put(2.3,4.3){\line(0,1){1}} \put(2.3,5.3){\line(0,1){1}}
\put(3.3,4.3){\line(0,1){6}} \put(4.3,4.3){\line(0,1){6}}
\put(5.3,4.3){\line(0,1){1}} \put(7.3,4.3){\line(0,1){6}}
\put(8.3,4.3){\line(0,1){6}} \put(3.3,5.3){\line(0,1){1}}
\put(4.3,5.3){\line(0,1){1}} \put(5.3,5.3){\line(0,1){1}}
\put(7.3,5.3){\line(0,1){1}} \put(8.3,5.3){\line(0,1){1}}
\put(2.3,8.3){\line(0,1){1}} \put(2.3,9.3){\line(0,1){1}}
\put(2.3,10.3){\line(0,1){0.5}}

\put(2.3,3.3){\line(1,-2){1.5}} \put(2.3,2.8){\line(1,-2){1.25}}
\put(2.3,2.3){\line(1,-2){1}} \put(2.3,1.8){\line(1,-2){0.75}}
\put(2.3,1.3){\line(1,-2){0.5}} \put(2.3,0.8){\line(1,-2){0.25}}
\put(2.55,3.3){\line(1,-2){1.5}} \put(2.8,3.3){\line(1,-2){1.5}}
\put(3.05,3.3){\line(1,-2){1.5}} \put(3.3,3.3){\line(1,-2){1.5}}
\put(3.55,3.3){\line(1,-2){0.75}} \put(3.8,3.3){\line(1,-2){0.5}}
\put(4.05,3.3){\line(1,-2){0.25}} \put(4.3,1.3){\line(1,-2){0.5}}
\put(4.55,1.3){\line(1,-2){0.5}} \put(4.8,1.3){\line(1,-2){0.5}}
\put(5.05,1.3){\line(1,-2){0.5}} \put(5.3,1.3){\line(1,-2){0.5}}
\put(5.55,1.3){\line(1,-2){0.5}} \put(5.8,1.3){\line(1,-2){0.5}}
\put(6.05,1.3){\line(1,-2){0.5}} \put(6.3,1.3){\line(1,-2){0.5}}
\put(6.55,1.3){\line(1,-2){0.5}} \put(2.3,1.30){\circle*{0.1}}
\put(3.3,1.30){\circle*{0.1}} \put(4.3,1.30){\circle*{0.1}}
\put(5.3,1.30){\circle*{0.1}} \put(2.3,2.3){\circle*{0.1}}
\put(3.3,2.3){\circle*{0.1}} \put(4.3,2.3){\circle*{0.1}}
\put(2.3,3.3){\circle*{0.1}} \put(3.3,3.3){\circle*{0.1}}
\put(4.3,3.3){\circle*{0.1}} \put(2.3,10.3){\circle*{0.1}}
\put(2.3,1.30){\circle*{0.1}} \put(3.3,1.30){\circle*{0.1}}
\put(4.3,1.30){\circle*{0.1}} \put(5.3,1.30){\circle*{0.1}}
\put(7.3,0.30){\circle*{0.1}} \put(2.3,0.30){\circle*{0.1}}
\put(3.3,0.30){\circle*{0.1}} \put(4.3,0.30){\circle*{0.1}}
\put(5.3,0.30){\circle*{0.1}} \put(2.05,2.05){\makebox(0,0){$2$}}
\put(2.05,1.05){\makebox(0,0){$1$}}
\put(5.15,0.00){\makebox(0,0){$3$}}
\put(4.15,0.00){\makebox(0,0){$2$}}
\put(3.15,0.00){\makebox(0,0){$1$}}
\put(2.1,0.1){\makebox(0,0){$0$}}
\put(8.1,0.00){\makebox(0,0){$(r/2\,, 0)$}}
\put(7.05,0.00){\makebox(0,0){$F$}}
\put(2.05,3.05){\makebox(0,0){$3$}}
\put(2.6,10.4){\makebox(0,0){$B$}}
\put(3.5,10.4){\makebox(0,0){$(0,r)$}}

\put(6.95,1.6){\makebox(0,0){$E$}}
\put(5.95,3.6){\makebox(0,0){$D$}}
\put(4.55,6.55){\makebox(0,0){$C$}}
\put(4.55,1.6){\makebox(0,0){$G$}} \put(2.5,3.6){\makebox(0,0){$A$}}
\put(4.55,3.65){\makebox(0,0){$H$}}
\end{picture}
\begin{center}
Fig. 1
\end{center}

Let us consider the following four triangles:  $P_{1}:=
\bigtriangleup\,BOF$,\, $P_{2}:= \bigtriangleup\,BAD$,\, $P_{3}:=
\bigtriangleup\,CGE$, \, and\, $P_{4}:= \bigtriangleup\,CHD$. It is
easy to see that $P_{1}$ is the $r$th dilate of the triangle

\bigskip

\bigskip

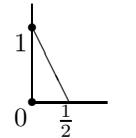
\setlength{\unitlength}{1cm}
\begin{picture}(3,0.5)

\put(8.15,0.1){\makebox(0,0){$0$}}
\put(8.15,1.1){\makebox(0,0){$1$}}
\put(8.75,0.05){\makebox(0,0){$\frac{1}{2}$}}
\put(9.75,0.30){\makebox(0,0){.}} \put(8.3,1.3){\line(1,-2){0.5}}
\put(8.3,0.3){\line(0,1){1.3}} \put(8.3,0.3){\line(1,0){1}}

\put(8.3,0.30){\circle*{0.1}} \put(8.3,1.30){\circle*{0.1}}
\end{picture}

\medskip

The direct computation shows that $L(P_{1}, 0) = 1$, $L(P_{1}, 1) =
2$, $L(P_{1}, 2) = 4$, and $L(P_{1}, 3) = 6$. As in Example 2, we
obtain that $$L(P_{1}, r) = {\frac{1}{4}}r^{2} + r + \left[1,
{\frac{3}{4}}\right]_{r}.$$

As it is seen from Fig. 1, the triangles $P_{2}$, $P_{3}$ and
$P_{4}$ are similar to $P_{1}$. If $N_{i}$ denotes the number of
lattice points in $P_{i}$ ($i=1, 2, 3, 4$), then for all
sufficiently large $r\in \mathbb{N}$,\,\,
$\chi^{(w)}_{\mathcal{A}}(r) = N_{1} - N_{2} - N_{3} + N_{4}$.
Therefore,
$$\chi^{(w)}_{\mathcal{A}}(r) = L(P_{1}, r) - L(P_{1}, r-3) - L(P_{1}, r-5) +
L(P_{1}, r-7)= {\frac{1}{2}}r + \left[5, {\frac{9}{2}}\right]_{r}.
$$  }
\end{example}

\section{The main result}

Let $K$ be a difference field of zero characteristic with a basic
set of translations $\sigma = \{\alpha_{1},\dots, \alpha_{m}\}$ that
are assigned positive integer weights $w_{1}, \dots, w_{m}$,
respectively. As before, let $T$ denote the free commutative
semigroup generated by the set $\sigma$. For any transform $\tau =
\alpha_{1}^{k_{1}}\dots\alpha_{m}^{k_{m}}\in T$, the number
$$\ord_{w}\tau = \D\sum_{i=1}^{m}w_{i}k_{i}$$ will be called the
{\em order} of $\tau$; it will be denoted by $\ord_{w}\tau$.

Furthermore, for any $r\in \mathbb{N}$, we set $$T_{w}(r) =
\{\tau\in T\,|\,\ord_{w}\tau\leq r\}.$$ The following theorem
establishes the existence of a dimension quasi-polynomial associated
with a difference field extension with weighted basic translations.

\begin{theorem}
With the above notation, let $L = K\langle
\eta_{1},\dots,\eta_{n}\rangle$ be a $\sigma$-field extension of $K$
generated by a finite set $\eta =\{\eta_{1},\dots,\eta_{n}\}$. For
any $r\in \mathbb{N}$, let $L_{r} =
K(\cup_{i=1}^{n}T_{w}(r)\eta_{i})$. Then there exists a
quasi-polynomial $\Phi^{(w)}_{\eta|K}(t)$ such that

\medskip

{\em (i)} \,$\Phi^{(w)}_{\eta|K}(r) = \trdeg_{K}L_{r}$ for all
sufficiently large $r\in \mathbb{N}$.

\medskip

{\em (ii)}\, $\deg \Phi^{(w)}_{\eta|K}\leq m = \Card\,\sigma$.

\medskip

{\em (iii)}\, $\Phi^{(w)}_{\eta|K}$ is an alternative sum of Ehrhart
quasi-polynomials associated with rational conic polytopes.

\medskip

{\em (iv)} \, The degree and leading coefficient of the
quasi-polynomial $\Phi^{(w)}_{\eta|K}$ are constants that do not
depend on the set of difference generators $\eta$ of the extension
$L/K$. Furthermore, the coefficient of $t^{m}$ in
$\Phi^{(w)}_{\eta|K}$ can be represented as $\D\frac{a}{m!w_{1}\dots
w_{m}}$ where $a$ is equal to the difference transcendence degree of
$L/K$.
\end{theorem}

In order to prove this theorem we need some results on reduction and
autoreduced sets in the ring of difference polynomials. A detailed
description of this technique, as well as proofs of Proposition 1
and Proposition 2 below,  can be found in \cite[Section 3.3]{KLMP}
and \cite[Section 2.4]{Levin2}.

Let $K$ be a difference field with a basic set $\sigma =
\{\alpha_{1},\dots, \alpha_{m}\}$, $T$ the free commutative
semigroup generated by $\sigma$, and  $R = K\{y_{1},\dots, y_{n}\}$
the algebra of difference ($\sigma$-) polynomials over $K$. As
before we set  $TY = \{\tau y_{i} | \tau \in T, 1\leq i\leq n\}$;
elements of this sets will be called {\em terms}. The order of a
term $\tau y_{i}$ (with respect to the given weights) is defined as
the order of $\tau$: $\ord_{w}(\tau y_{i}) = \ord_{w}\tau$.

By a {\em ranking} of the set $TY$ we mean a well-ordering $\leq$ of
this set satisfying the following two conditions.

\smallskip

(i) $u\leq \tau u$ for any $u\in TY, \tau \in T$.

\smallskip

(ii) If $u, v\in TY$ and $u\leq v$, then $\tau u\leq \tau v$ for
 any $\tau \in T$.

\smallskip

In the rest of the paper we consider a ranking defined as follows:

If $\tau = \alpha_{1}^{k_{1}}\dots\alpha_{m}^{k_{m}},\, \tau' =
\alpha_{1}^{l_{1}}\dots\alpha_{m}^{l_{m}}\in T$ and $i,
j\in\{1,\dots, n\}$, then $\tau y_{i} < \tau'y_{j}$ if and only if
 $$(\ord_{w}\tau, k_{1},\dots, k_{m}, i) <_{lex} (\ord_{w}\tau', l_{1},\dots, l_{m},
 j).$$

If $A\in K\{y_{1},\dots, y_{n}\}\setminus K$, then the greatest
(with respect to our ranking) element of $TY$ that appears in $A$ is
called the {\em leader} of $A$; it is denoted by $u_{A}$. If $A$ is
written as a polynomial in $u_{A}$, $A =
\sum_{i=0}^{d}I_{i}u_{A}^{i}$ ($d=deg_{u_{A}}A$ and the
$\sigma$-polynomials $I_{0},\dots, I_{d}$ do not contain $u_{A}$),
then $I_{d}$ is called the {\em initial} of the $\sigma$-polynomial
$A$; it is denoted by $I_{A}$.

Let $A$, $B\in K\{y_{1},\dots, y_{n}\}$. We say that $A$ has lower
rank than $B$ and write $A < B$, if either $A\in K, B\notin K$ or
$u_{A}< u_{B}$ or $u_{A}= u_{B} = u$, $deg_{u}A < deg_{u}B$. If
neither $A < B$ nor $B < A$, we say that $A$ and $B$ have the same
rank and write $rk\,A = rk\,B$. Furthermore, we say that $A$ is {\em
reduced} with respect to $B$ ($B\notin K$) if $A$ does not contain
any power of a transform $\tau u_{B}$ ($\tau \in T$) whose exponent
is greater than or equal to $deg_{u_{B}}B$. If $S$ is any subset of
$K\{y_{1},\dots, y_{n}\}\setminus K$, then a $\sigma$-polynomial
$A\in K\{y_{1},\dots, y_{n}\}$ is said to be reduced with respect to
$S$ if $A$ is reduced with respect to every element of $S$.

\smallskip

A set $\Sigma \subseteq K\{y_{1},\dots, y_{n}\}$ is called an {\em
autoreduced set} if either $\Sigma = \emptyset$ or $\Sigma \bigcap K
= \emptyset$  and every element of $\Sigma$ is reduced with respect
to all other elements of $\Sigma$. As it is shown in\cite[Section
2.4]{Levin2}, every autoreduced set is finite and distinct elements
of an autoreduced set have distinct leaders.

\begin{proposition} Let  $\cal{A}$ $ = \{A_{1},\dots,
A_{p}\}$ be an autoreduced  set in $K\{y_{1},\dots, y_{n}\}$. Let
$I({\cal{A}}) = \{B\in K\{y_{1},\dots, y_{n}\}  |$ either $B =1$ or
$B$ is a product of finitely many $\sigma$-polynomials of the form
$\tau(I_{A_{i}})$ ($\tau \in T, i=1,\dots, p)\}$. Then for any $C\in
K\{y_{1},\dots, y_{n}\}$, there exist $\sigma$-polynomials $J\in
I(\cal{A})$  and $C_{0}$ such that $C_{0}$ is  reduced  with respect
to the set $\cal{A}$ and $JC\equiv C_{0}(mod\,[\cal{A}])$, that is,
$JC-C_{0}\in [\cal{A}]$.
\end{proposition}

Let $\cal{A}$ $ = \{A_{1},\dots, A_{p}\}$ and $\cal{B}$ $ =
\{B_{1},\dots, B_{q}\}$ be two autoreduced sets whose elements are
written in the order of increasing rank. We say that $\cal{A}$ has
lower rank than $\cal{B}$ if one of the following conditions holds:

(i) {\em there exists $k\in\mathbb{N},\,  1\leq k\leq \min\{p, q\}$,
such that $rk\,A_{i} = rk\,B_{i}$ for $i=1,\dots, k-1$ and $A_{k} <
B_{k}$;}

(ii) \,{\em  $p > q$ and  $rk\,A_{i} = rk\,B_{i}$ for $i=1,\dots,
q$.}

\bigskip

\begin{proposition} Every nonempty set of autoreduced sets
contains an autoreduced set of lowest rank. (If $J$ is a difference
ideal of $K\{y_{1},\dots, y_{n}\}$, an autoreduced subset of $J$ of
lowest rank is called a {\bf characteristic set} of $J$.)

Let $\Sigma$ be a characteristic set of a $\sigma$-ideal $J$ of
$K\{y_{1},\dots, y_{n}\}$. Then $J$ does not contain nonzero
difference polynomials reduced with respect to $\Sigma$. In
particular, if $A\in \Sigma$, then $I_{A}\notin J$.
\end{proposition}

\medskip

{\centerline{PROOF OF THEOREM 4}}

\bigskip

Let $P$ be the defining $\sigma$-ideal of the extension $L =
K\langle\eta_{1},\dots, \eta_{n}\rangle$, that is $P =
\Ker(K\{y_{1},\dots, y_{n}\}\rightarrow L),\,\,\,\, y_{i}\mapsto
\eta_{i}.$

Let $\Sigma = \{A_{1},\dots, A_{d}\}$ be a characteristic set of
$P$, let $u_{i}$ denote the leader of $A_{i}$ ($1\leq i\leq d$) and
for every $j=1,\dots, n$, let $$E_{j} = \{(k_{1},\dots,
k_{m})\in\mathbb{N}^{m} | \alpha_{1}^{k_{1}}\dots
\alpha_{m}^{k_{m}}y_{j}\,\,\, \text{is a leader of a}\,\,\,
\sigma\text{-polynomial in}\,\,\,  {\cal{A}}\}.$$

Let $V = \{u\in TY\,|\,u$ is not a transform of any $u_{i}$ ($1\leq
i\leq s$) and for every $r\in\mathbb{N}$, let $V(r) = \{u\in
V\,|\,\ord_{w}\,u\leq r\}$. By Proposition 2, the ideal $P$ does not
contain non-zero difference polynomials reduced with respect to
${\cal{A}}$. It follows that for every $r\in\mathbb{N}$, the set
$V_{\eta}(r) = \{v(\eta)\,|\,v\in V(r)\}$ is algebraically
independent over $K$. Indeed, if there exists a nonzero polynomial
$B\in K[X_{1},\dots, X_{k}]$ ($k> 0$) in $k$ variables over $K$ and
elements $v_{1},\dots, v_{k}\in V_{\eta}(r)$ such that
$B(v_{1}(\eta),\dots, v_{k}(\eta)) = 0$. Then $B(v_{1},\dots,
v_{k})\in P$ and this $\sigma$-polynomial is reduced with respect to
${\cal{A}}$, so we arrive at a contradiction.

If $A_{i}\in\Sigma$, then $A_{i}(\eta) = 0$, hence $u_{i}(\eta)$ is
algebraic over the field $K(\{\tau\eta_{j}\,|\,\tau y_{j} < u_{i}\,
(\tau \in T, \,1\leq j\leq n)\})$. Therefore, for any $\theta\in T$,
$\theta u_{i}$ is algebraic over the field
$K(\{\tau\eta_{j}\,|\,\tau y_{j} < \theta u_{i}\, (\tau \in T,
\,1\leq j\leq n)\})$. By induction on the well-ordered set $TY$ (and
using the fact that if $u, v\in TY$ and $u < v$, then $\ord_{w}\leq
\ord_{w}v$), we obtain that for every $r\in\mathbb{N}$, the field
$L_{r} = K(\{\tau\eta_{j}\,|\,\tau\in T(r),\, 1\leq j\leq n\})$ is
an algebraic extension of the field $K(\{v(\eta)\,|\,v\in V(r)\})$.
It follows that $V_{\eta}(r)$ is a transcendence basis of $L_{r}$
over $K$ and $\trdeg_{K}L_{r} = \Card\,V_{\eta}(r)$.

\smallskip

The number of terms $\alpha_{1}^{k_{1}}\dots \alpha_{m}^{k_{m}}
y_{j}$ in $V(r)$ is equal to the number of $m$-tuples $k =
(k_{1},\dots, k_{m})\in\mathbb{N}^{m}$ such that $\ord_{w}k\leq r$
and $k$ does not exceed any $m$-tuple in $E_{j}$ with respect to the
product order on $\mathbb{N}^{m}$. By Theorem 3, this number is
expressed by a quasi-polynomial of degree at most $m$.

Therefore, for all sufficiently large $r\in \mathbb{N}$,
$$\trdeg_{K}L_{r} = \Card\,V_{\eta}(r) =
\D\sum_{j=1}^{n}\chi^{(w)}_{E_{j}}(r)$$ where
$\chi^{(w)}_{E_{j}}(t)$ is the dimension quasi-polynomial of the set
$E_{j}\subseteq\mathbb{N}^{m}$. Since each $\chi^{(w)}_{E_{j}}(t)$
is an alternative sum of Ehrhart quasi-polynomials associated with
conic polytopes (see formula (1) in Theorem 3),
$\Phi^{(w)}_{\eta|K}(t) = \D\sum_{j=1}^{n}\chi^{(w)}_{E_{j}}(t)$ has
a similar representation and satisfies conditions (i) -- (iii) of
Theorem 4.

\medskip

If $\zeta = (\zeta_{1},\dots, \zeta_{k})$ is another system of
$\sigma$-generators of the extension $L/K$, so that $L =
K\langle\eta_{1},\dots, \eta_{n}\rangle = K\langle\zeta_{1},\dots,
\zeta_{k}\rangle$, then there exists $q\in \mathbb{N}$ such that
$\eta_{1},\dots, \eta_{n}\in K(\cup_{i=1}^{k}T_{w}(q)\zeta_{i})$ and
$\zeta_{1},\dots, \zeta_{k}\in K(\cup_{i=1}^{n}T_{w}(q)\eta_{i})$.
Therefore, for all sufficiently large $r\in \mathbb{N}$ (namely, for
all $r\geq q$), one has
$$K(\cup_{i=1}^{n}T_{w}(r)\eta_{i})\subseteq
K(\cup_{i=1}^{k}T_{w}(r+q)\zeta_{i})\,\,\, \text{and}\,\,\,
K(\cup_{i=1}^{k}T_{w}(r)\zeta_{i})\subseteq
K(\cup_{i=1}^{n}T_{w}(r+q)\eta_{i}),$$ that is,
$\Phi^{(w)}_{\eta|K}(r)\leq \Phi^{(w)}_{\zeta|K}(r+q)$ and
$\Phi^{(w)}_{\zeta|K}(r)\leq \Phi^{(w)}_{\eta|K}(r+q)$. It follows
that the quasi-polynomials $\Phi^{(w)}_{\eta|K}(t)$ and
$\Phi^{(w)}_{\zeta|K}(t)$ have equal degrees and equals leading
coefficients.

\medskip

In order to prove the last statement of the theorem, note first that
if the elements $\eta_{1},\dots,\eta_{n}$ are $\sigma$-algebraically
independent over $K$, then $$\Phi^{(w)}_{\eta|K}(t) =
n\sum\lambda^{(m)}_{w}(t).$$

Indeed, if $r\in \mathbb{N}$, and $\Xi_{i}(r) = \{\xi =
\alpha_{1}^{k_{1}}\dots\alpha_{m}^{k_{m}}\eta_{i}\,|\,k_{i}\in
\mathbb{N},\,\,\ord_{w}\xi\leq r\}$ for $i = 1,\dots, n$, then
$\cup_{i=1}^{n}\Xi_{i}(r)$ is a transcendence basis of the field
extension $K(\cup_{i=1}^{n}T_{w}(r)\eta_{i})/K$ and the number of
elements of this basis is equal to $n\Card\{(k_{1},\dots,
k_{m})\in\mathbb{N}^{m}\,|\,\sum_{j=1}^{m}w_{j}k_{j}\leq r\} =
n\sum\lambda^{(m)}_{w}(r)$.

\medskip

Let $L = K\langle\eta_{1},\dots,\eta_{n}\rangle$ and $d =
\trdeg_{K}L$. Without loss of generality we can assume that
$\eta_{1},\dots, \eta_{d}$ is a $\sigma$-transcendence basis of
$L/K$. (We allow $d=0$; in this case every element of $L$ is
$\sigma$-algebraic over $K$.) Then for any $j=d+1,\dots, n$, the
element $\eta_{j}$ is $\sigma$-algebraic over the $\sigma$-field $F
= K\langle\eta_{1},\dots,\eta_{d}\rangle$ (properties of
$\sigma$-transcendence dependence and $\sigma$-transcendence bases,
as it is shown in \cite[Section 4.1]{Levin2}, are similar to
properties of usual transcendence dependence and transcendence bases
of field extensions).

Let $B_{j}\in F\{y_{j}\}$ be a nonzero $\sigma$-polynomial of the
smallest total degree such that $B_{j}(\eta_{j}) = 0$. If $u_{j} =
\tau_{j}y_{j}$ is the leader of $B_{j}$ ($1\leq k_{j}\leq n$) and
the $\sigma$-polynomial $B_{j}$ is written as a polynomial in
$u_{j}$, $B_{j} = \sum_{i=1}^{k}I_{ij}u_{j}^{i}$, then
\begin{equation}
B_{j}(\eta_{j}) = \sum_{i=1}^{k}I_{ij}(\eta_{j})u_{j}^{i}(\eta_{j})
= 0.
\end{equation}
By the choice of $B_{j}$, $I_{kj}(\eta_{j})\neq 0$, so that
$u_{j}(\eta_{j})$ is algebraic over the field
$F(\{\tau\eta_{j}\,|\,\tau\in T$ and $\tau y_{j} < u_{j}\})$. (Here
$<$ denotes the considered total order on the set of terms $TY$).
Therefore, we can fix $p_{j}\in\mathbb{N}$ such that $u(\eta_{j})$
is algebraic over the field $K(\{\theta\eta_{i}\,|\,1\leq i\leq d$
and $\theta\in T_{w}(p_{j})\} \cup \{\tau\eta_{j}\,|\,\tau\in
T$\,and\, $\tau y_{j} < u_{j}\})$.

\smallskip

It follows from equality (2) that for any $\tau\in T$,
$\sum_{i=1}^{k}(\tau I_{ij}(\eta_{j}))(\tau u_{j}(\eta_{j}))^{i} =
0$, so that $\tau u_{j}(\eta_{j})$ is algebraic over the field
$K(\{\theta\eta_{i}\,|\,\theta\in T_{w}(p_{j}+\ord_{w}\tau),\, 1\leq
i\leq d\}\cup\{\theta'\eta_{j}\,|\,\theta'\in T,\, \theta' y_{j} <
\tau u_{j}\}$. By induction on the well-ordered set $TY$ we obtain
that $\tau u_{j}(\eta_{j})$ is algebraic over the field
$K(\{\theta\eta_{i}\,|\,\theta\in T_{w}(p_{j}+\ord_{w}\tau),\, 1\leq
i\leq d\}\cup\{\theta'\eta_{j}\,|\,\theta'\in T, \theta'y_{j} < \tau
u_{j}$ and $\tau_{j}\nmid\theta'\})$.

Let $p = \max\{p_{d+1},\dots, p_{n}\}$ and
$s\geq\max\{\ord_{w}u_{j}\,|\,d+1\leq j\leq n\}$. Then every element
of the field $$L_{s} = K(\{\tau\eta_{i}\,|\,\tau\in
T_{w}(s)\eta_{i},\, 1\leq i\leq n\})$$ is algebraic over the field

\medskip

\noindent$M_{s} = K(\{\tau\eta_{i}\,|\,\tau\in T_{w}(s+p), 1\leq
i\leq d\}\cup\{\theta\eta_{l}\,|\,\theta\in[T_{w}(s) -
T_{w}(s-r_{l})],$

\noindent$d+1\leq l\leq n\})$ where $r_{l} = \ord_{w}u_{l}$
($d+1\leq l\leq n$). It follows that

\medskip

\noindent$\trdeg_{K}L_{s}\leq\trdeg_{K}L_{s}M_{s} = \trdeg_{K}M_{s}
+ \trdeg_{M_{s}}L_{s}M_{s} =\trdeg_{K}M_{s}.$

\medskip

Since $\trdeg_{K}L_{s} = \Phi^{(w)}_{\eta|K}(s)$ and
$$\trdeg_{K}M_{s}\leq d\Card\,T_{w}(s+p)
+\sum_{l=d+1}^{n}\Card\,[T_{w}(s) - T_{w}(s-r_{l})],$$ we have
$$\Phi^{(w)}_{\eta|K}(s)\leq d\lambda_{w}^{(m)}(s+p) +
\sum_{l=d+1}^{n}[\lambda_{w}^{(m)}(s) -
\lambda_{w}^{(m)}(s-r_{l})].$$ Since $\deg(\lambda_{w}^{(m)}(t) -
\lambda_{w}^{(m)}(t-r_{j})) < m$, we obtain that the coefficient of
$t^{m}$ in $\Phi^{(w)}_{\eta|K}(t)$ does not exceed
$\D\frac{d}{m!w_{1}\dots w_{m}}$. On the other hand, the set $\{\tau
\eta_{i}\,|\,\tau\in T_{w}(s),\, 1\leq i\leq d\}$ is
$\sigma$-algebraically independent over $K$, hence
$\Phi^{(w)}_{\eta|K}(s)\geq d\Card\,T_{w}(s) =
d\lambda_{w}^{(m)}(s)$, so the coefficient of $t^{m}$ in
$\Phi^{(w)}_{\eta|K}(t)$ is greater than or equal to
$\D\frac{d}{m!w_{1}\dots w_{m}}$. Thus, if the coefficient of
$t^{m}$ in $\Phi^{(w)}_{\eta|K}$ is represented as
$\D\frac{a}{m!w_{1}\dots w_{m}}$, then $a = d =
\sigma$-$\trdeg_{K}L$.\, $\square$

\bigskip

Theorem 4 allows one to assign a quasi-polynomial to a system of
algebraic difference equations with weighted basic translations
\begin{equation}
f_{i}(y_{1},\dots, y_{n}) = 0\hspace{0.5in}(i=1,\dots, p)
\end{equation}
($f_{i},\in R = K\{y_{1},\dots, y_{n}\}$ for $i=1,\dots, p$) such
that the difference ideal $P$ of $R$ generated by the
$\sigma$-polynomials $f_{1}, \dots, f_{p}$ is prime (e. g. to a
system of linear difference equations). Systems of this form arise,
in particular, as finite difference approximations of systems of
PDEs with weighted derivatives (see, for example, \cite{Shananin1}
and \cite{Shananin2}).

\medskip

In this case, one can consider the quotient field $L = \qf(R/P)$ as
a finitely generated $\sigma$-field extension of $K$: $L =
K\langle\eta_{1},\dots, \eta_{n}\rangle$ where $\eta_{i}$ is the
canonical image of $y_{i}$ in $R/P$. The corresponding dimension
quasi-polynomial $\Phi^{(w)}(t) = \Phi^{(w)}_{\eta|K}(t)$ is said to
be the {\em difference dimension quasi-polynomial} of system (3). It
can be viewed as the Einstein's strength of the system of partial
difference equations with weighted translations (see \cite[Section
7.7]{Levin2} for a detail description of this concept).

\begin{example}
{\em Let $K$ be a difference field with basic set of translations
$\sigma = \{\alpha_{1}, \alpha_{2}\}$ where $\alpha_{1}$ and
$\alpha_{2}$ are assigned weights $2$ and $1$, respectively.
Consider any system of linear algebraic difference equations in one
difference indeterminate $y$ over $K$ that consist of two equations
\begin{equation}
f_{i}(y) = 0\hspace{0.5in}  (i=1, 2)
\end{equation}
such that the leaders of the $\sigma$-polynomials $f_{1}$ and
$f_{2}$ are $\alpha_{1}^{2}\alpha_{2}y$ and $\alpha_{2}^{3}y$,
respectively. Since $[f_{1}, f_{2}]$ is a linear difference ideal of
$K\{y\}$, it is prime (see \cite[Proposition 2.4.9]{Levin2}) and
$\{f_{1}, f_{2}\}$ is its characteristic set. Using the result of
Example 2, we obtain the following expression for the difference
dimension quasi-polynomial of system (4):
$$\Phi^{(w)}(t) = {\frac{1}{2}}t + \left[5,
{\frac{9}{2}}\right]_{t}.$$ }
\end{example}

\end{document}